\def\Bbb#1{{\fam\msbfam\relax#1}}
\font\fivemsb=msbm5
\font\sevenmsb=msbm7
\font\tenmsb=msbm10
\newfam\msbfam
\textfont\msbfam=\tenmsb
\scriptfont\msbfam\sevenmsb
\scriptscriptfont\msbfam\fivemsb

\def\spc{{\Bbb C}}
\def\spr{{\Bbb R}}
\def\vf{\varphi}
\def\wt{\widetilde}
\def\wh{\widehat}
\def\a{\alpha}
\def\b{\beta}
\def\d{\delta}

\def\F{\Phi}

\def\g{\gamma}
\def\l{\lambda}
\def\m{\mu}
\def\r{\rho}
\def\s{\sigma}
\def\p{\pi}

\def\t{\theta}
\def\x{\xi}
\def\y{\eta}

\def\ce{{\cal E}}

\def\cl{{\cal L}}

\def\la{\langle}
\def\ra{\rangle}
\def\lla{\langle\!\langle}
\def\rra{\rangle\!\rangle}

\documentstyle[verbatim]{amsart}

\theoremstyle{plain}
\newtheorem{theorem}{Theorem}[section]

\newtheorem{lemma}[theorem]{Lemma}

\theoremstyle{definition}

\newtheorem{example}[theorem]{Example}

\numberwithin{equation}{section}

\begin{document}
\title[CKS-space, growth functions]%
   {CKS-space in terms of growth functions}

\author{Nobuhiro Asai}
\author{Izumi Kubo}
\author{Hui-Hsiung Kuo}
\address{Nobuhiro Asai: Graduate School of Mathematics\\
  Nagoya University \\ Nagoya 464-8602 \\ JAPAN}
\address{Izumi Kubo: Department of Mathematics \\
Graduate School of Science \\ Hiroshima University \\
Higashi-Hiroshima 739-8526 \\ JAPAN}
\address{Hui-Hsiung Kuo: Department of Mathematics\\
  Louisiana State University \\ Baton Rouge \\
LA 70803 \\ USA}

\maketitle

\begin{abstract}
A class of growth functions $u$ is introduced to construct
Hida distributions and test functions. The Legendre
transform $\ell_{u}$ of $u$ is used to define a sequence
$\a(n)=(\ell_{u}(n) n!)^{-1},\,n\geq 0$, of positive
numbers. From this sequence we get a CKS-space. Under
various conditions on $u$ we show that the associated
sequence $\{\a(n)\}$ satisfies those conditions for
carrying out the white noise distribution theory on the
CKS-space. We show that $u$ and its dual Legendre transform
$u^{*}$ are growth functions for test and generalized
functions, respectively, in the characterization theorems.

\end{abstract}

\section{CKS-space} \label{sec:1}

Let $\{\a(n)\}_{n=0}^{\infty}$ be a sequence of positive
numbers. The exponential generating functions $G_{\a}$ and
$G_{1/\a}$ of the sequences $\{\a(n)\}$ and $\{1/\a(n)\}$
are defined by
\begin{equation} \label{eq:1-1}
G_{\a}(r) = \sum_{n=0}^{\infty} {\a(n) \over n!} r^{n},
   \qquad G_{1/\a}(r) = \sum_{n=0}^{\infty} {1 \over
    n!\a(n)} r^{n}.
\end{equation}
Later on we will impose various conditions from the following
list:

\medskip
\begin{itemize}
\item[(A1)] $\a(0)=1$ and $\inf_{n\geq 0} \a(n) \s^{n}
>0$ for some $\s\geq 1$.
\medskip
\item[(A2)] $\lim_{n\to\infty} \left({\a(n) \over n!}
\right)^{1/n} =0$.
\medskip
\item[($\wt{\text{A}}2$)] $\lim_{n\to\infty} \left({1\over
  n!\a(n)}\right)^{1/n} =0$.
\smallskip
\item[(B1)] $\limsup_{n\to\infty} \left({n! \over \a(n)}
\inf_{r>0} {G_{\a}(r) \over r^{n}}\right)^{1/n} < \infty$.
\smallskip
\item[($\wt{\text{B}}1$)] $\limsup_{n\to\infty} \left(
n!\a(n) \inf_{r>0} {G_{1/\a}(r) \over r^{n}}\right)^{1/n}
  < \infty$.
\medskip
\item[(B2)] The sequence $\g(n)={\a(n) \over n!}, n\geq 0,$
is log-concave, i.e., for all $n\geq 0$,
\begin{equation}
\g(n) \g(n+2) \leq \g(n+1)^{2}. \notag
\end{equation}
\item[($\wt{\text{B}}2$)] The sequence $\left\{{1 \over
n!\a(n)}\right\}$ is log-concave.
\bigskip
\item[(B3)] The sequence $\{\a(n)\}$ is log-convex, i.e.,
for all $n\geq 0$,
\begin{equation}
\a(n) \a(n+2) \geq \a(n+1)^{2}. \notag
\end{equation}
\item[(C1)] There exists a constant $c_{1}$ such that for
all $n\leq m$,
\begin{equation}
\a(n) \leq c_{1}^{m} \a(m).  \notag
\end{equation}
\item[(C2)] There exists a constant $c_{2}$ such that for
all $n$ and $m$,
\begin{equation}
\a(n+m) \leq c_{2}^{n+m} \a(n) \a(m).  \notag
\end{equation}
\item[(C3)] There exists a constant $c_{3}$ such that for
all $n$ and $m$,
\begin{equation}
\a(n) \a(m) \leq c_{3}^{n+m} \a(n+m).  \notag
\end{equation}
\end{itemize}

Recently Cochran et al.~\cite{cks} have used a sequence
$\{\a(n)\}$ of positive numbers to define spaces of test
and generalized functions on a white noise space. They
assumed condition (A1) with $\s=1$. But our (A1) is strong
enough to imply that the space of test functions is
contained in the $L^{2}$-space of the white noise measure.
In \cite{cks} conditions (A2) (B1) (B2) are considered.
Condition (A2) is to assure that the function $G_{\a}$
in Equation (\ref{eq:1-1}) is an entire function.
Condition (B1) is used for the characterization theorem of
generalized functions in Theorems 5.1 and 6.1 \cite{cks}.
Condition (B2) is shown to imply condition (B1) in
Theorem 4.3 \cite{cks}.

In the papers by Asai et al.~\cite{akk1} \cite{akk2},
conditions ($\wt{\text{A}}2$) ($\wt{\text{B}}1$)
($\wt{\text{B}}2$) (B3) are considered. It can be
easily checked that condition (A1) implies condition
($\wt{\text{A}}2$). Condition ($\wt{\text{A}}2$) is to
assure that the function $G_{1/\a}$ in Equation
(\ref{eq:1-1}) is an entire function. In \cite{akk2}
condition ($\wt{\text{B}}1$) is used for the
characterization theorem of test functions. Condition
($\wt{\text{B}}2$) implies condition ($\wt{\text{B}}1$),
while obviously condition (B3) implies condition
($\wt{\text{B}}2$).

In the paper by Kubo et al.~\cite{kks}, conditions (C1)
(C2) (C3) are assumed in order to carry out the
distribution theory for a CKS-space. As pointed out in
\cite{kks}, condition (C3) implies condition (C1).

An important example of $\{\a(n)\}$ is the sequence
$\{b_{k}(n)\}$ of Bell's numbers of order $k\geq 2$. The
sequence $\{b_{k}(n)\}$ satisfies conditions (A1) (A2)
(B1) (as shown in \cite{cks}), (B2) (B3) (as shown in
\cite{akk1}) (C1) (C2) (C3) (as shown in \cite{kks}.)
Therefore, Bell's numbers satisfy all conditions in
the above list.

The essential conditions for distribution theory on a
CKS-space are (A1) (A2) (B2) ($\wt{\text{B}}2$) (C2)
(C3). All other conditions can be derived from these
six conditions except for (B3). We have taken
($\wt{\text{B}}2$) instead of (B3) for the following
reason. The condition (B3) is rather strong and we do
not know how to prove this condition for a growth
function $u$ in Section \ref{sec:3}. Fortunately, we do
not need (B3) for white noise distribution theory.

Now, we briefly explain the CKS-space associated with
a sequence $\{\a(n)\}$ of positive numbers. Let $\ce$ be
a nuclear space with topology given by a sequence of
inner product norms $\{|\cdot|_{p}\}_{p=0}^{\infty}$.
Let $\ce_{p}$ be the completion of $\ce$ with respect to
the norm $|\cdot|_{p}$. Assume the conditions:
\begin{itemize}
\item[(a)] There exists a constant $0<\r<1$ such that
$|\cdot|\leq \r|\cdot|_{1}\leq \cdots\leq \r^{p}|\cdot|_{p}
\leq \cdots$.
\item[(b)] For any $p\geq 0$, there exists some $q\geq p$
such that the inclusion mapping $i_{q, p}: \ce_{q} \to
\ce_{p}$ is a Hilbert-Schmidt operator.
\end{itemize}

\smallskip
Let $\ce'$ be the dual space of $\ce$. Then we have a
Gel'fand triple $\ce \subset \ce_{0} \subset \ce'$. Let
$\m$ be the standard Gaussian measure on $\ce'$ and let
$(L^{2})$ denote the Hilbert space of complex-valued
square integrable functions on $(\ce', \m)$. Every
$\vf\in (L^{2})$ can be uniquely represented by a sum of
multiple Wiener integrals
\begin{equation}
\vf=\sum_{n=0}^{\infty} I_{n}(f_{n}), \quad
f_{n} \in \ce_{0}^{\wh\otimes n}. \notag
\end{equation}
Moreover, the $(L^{2})$-norm $\|\vf\|_{0}$ of $\vf$ is
given as follows:
\begin{equation}
\|\vf\|_{0}^{2} = \sum_{n=0}^{\infty} n!|f_{n}|_{0}^{2}.
  \notag
\end{equation}

For a positive integer $p$ such that $\r^{-2p}\geq \s$, or
equivalently $p\geq (-2\log \r)^{-1} \log \s$, define
$\|\cdot\|_{p, \a}$ by
\begin{equation} \label{eq:1-2}
\|\vf\|_{p, \a}^{2} = \sum_{n=0}^{\infty} n! \a(n)
|f_{n}|_{p}^{2}.
\end{equation}
We can use the above assumption (a) to show that
\begin{equation}
\|\vf\|_{p, \a}^{2} \geq \left(\inf_{n\geq 0} \a(n) \s^{n}
\right) \|\vf\|_{0}^{2}.  \notag
\end{equation}
Thus if we define $[\ce_{p}]_{\a} = \{\vf\in (L^{2}); \,
\|\vf\|_{p, \a} < \infty\}$, then condition (A1) implies
that $[\ce_{p}]_{\a} \subset (L^{2})$ for all $p\geq
(-2\log \r)^{-1} \log \s$. The space $[\ce]_{\a}$ of test
functions on $\ce'$ is defined to be the projective limit
of $\{[\ce_{p}]_{\a}; \,p\geq (-2\log \r)^{-1} \log \s\}$.
Its dual space $[\ce]_{\a}^{*}$ is the space of generalized
functions on $\ce'$. Then we get the following Gel'fand
triple
\begin{equation} \label{eq:b}
[\ce]_{\a} \subset (L^{2}) \subset [\ce]_{\a}^{*}.
\end{equation}
This Gel'fand triple was introduced by Cochran et al.~in
\cite{cks} and is often referred to as the CKS-space
associated with the sequence $\{\a(n)\}$. Note that
$[\ce]_{\a}^{*} = \cup_{p} [\ce_{p}]_{\a}^{*}$ and the
$[\ce_{p}]_{\a}^{*}$-norm is given by
\begin{equation} \label{eq:1-3}
\|\vf\|_{-p, 1/\a}^{2} = \sum_{n=0}^{\infty} {n! \over
\a(n)} |f_{n}|_{-p}^{2}.
\end{equation}
For each $\x$ belonging to the complexification $\ce_{c}$
of $\ce$, define
\begin{equation}
:\!e^{\la \cdot, \x\ra}\!: \,= \sum_{n=0}^{\infty}
{1\over n!} I_{n}(\x^{\otimes n}).  \notag
\end{equation}
By Equations (\ref{eq:1-2}) and (\ref{eq:1-3}), we have
\begin{equation} \label{eq:1-4}
\|\!:\!e^{\la \cdot, \x\ra}\!:\!\|_{p, \a}^{2} =
G_{\a}\big(|\x|_{p}^{2}\big), \quad
\|\!:\!e^{\la \cdot, \x\ra}\!:\!\|_{-p, 1/\a}^{2} =
G_{1/\a}\big(|\x|_{-p}^{2}\big).
\end{equation}
The first equality shows that $:\!e^{\la \cdot, \x\ra}\!:
\,\in [\ce]_{\a}$ for all $\x\in\ce_{c}$.

The $S$-{\it transform} of a generalized function $\F$ in
$[\ce]_{\a}^{*}$ is defined to be the function
\begin{equation}
(S\F)(\x) = \lla \F, \, :\!e^{\la \cdot, \x\ra}\!:\rra,
\qquad \x \in \ce_{c},
\end{equation}
where $\lla \cdot, \cdot \rra$ denotes the bilinear pairing
of $[\ce]_{\a}^{*}$ and $[\ce]_{\a}$.

Under the condition (B1) it is proved in \cite{cks} that a
complex-valued function $F$ on $\ce_{c}$ is the
$S$-transform of a generalized function in $[\ce]_{\a}^{*}$
if and only if it satisfies an analyticity condition and
the growth condition: There exist constants
$K, a, p\geq 0$ such that
\begin{equation} \label{eq:1-5}
|F(\x)| \leq K G_{\a}\big(a|\x|_{p}^{2}\big)^{1/2},
  \qquad  \x\in\ce_{c}.
\end{equation}
Here the inequality is motivated by the first equality in
Equation (\ref{eq:1-4}).

On the other hand, under the condition ($\wt{\text{B}}1$),
it is shown in \cite{akk2} that a complex-valued function
$F$ on $\ce_{c}$ is the $S$-transform of a test function in
$[\ce]_{\a}$ if and only if it satisfies the analyticity
condition and the growth condition: For any $a, p\geq 0$,
there exists a constant $K\geq 0$ such that
\begin{equation} \label{eq:1-6}
|F(\x)| \leq K G_{1/\a}\big(a|\x|_{-p}^{2}\big)^{1/2},
  \qquad  \x\in\ce_{c}.
\end{equation}
This inequality is motivated by the second equality in
Equation (\ref{eq:1-4}).

Observe that the functions $G_{\a}$ and $G_{1/\a}$ are
defined by series as in Equation (\ref{eq:1-1}). In general
it is impossible to find the sums of these series in closed
forms. It is so even for the Kondratiev-Streit space
\cite{ks92} \cite{ks93} \cite{kuo96} when $\a(n)=(n!)^{\b},
\,\b\not= 0$. Thus it is desirable to find elementary
functions which can be used as growth functions in
Equations (\ref{eq:1-5}) and (\ref{eq:1-6}). This leads to
the investigation by Asai et al.~in the papers \cite{akk3}
\cite{akk4}. We will explain the essential ideas in the
rest of this paper.

\section{Legendre and dual Legendre transforms}
   \label{sec:2}

Let $C_{+, \log}$ denote the set of all positive continuous
functions $u$ on $[0, \infty)$ satisfying the condition:
\begin{equation}
\lim_{r\to\infty} {\log u(r) \over \log r}=\infty. \notag
\end{equation}

The {\it Legendre transform} $\ell_{u}$ of $u \in
C_{+, \log}$ is defined to be the function
\begin{equation}
\ell_{u}(t) = \inf_{r>0} {u(r) \over r^{t}}, \qquad
  t\in [0, \infty). \notag
\end{equation}

Below we state some properties of the Legendre transform.
For the proofs, see the paper \cite{akk3}.

\begin{theorem} \label{thm:2-1}
Let $u\in C_{+, \log}$. Then
\par\noindent
(1) The function $\ell_{u}$ is log-concave, i.e., for any
$r, s\geq 0$ and $0\leq \l\leq 1$,
\begin{equation}
\ell_{u}(\l r + (1-\l)s) \geq \ell_{u}(r)^{\l}
  \ell_{u}(s)^{1-\l}.  \notag
\end{equation}
In particular, the sequence $\{\ell_{u}(n)\}_{n=0}^{\infty}$
is log-concave.
\par\noindent
(2) For any nonnegative integers $n$ and $m$,
\begin{equation}
  \ell_{u}(0) \ell_{u}(n+m) \leq \ell_{u}(n) \ell_{u}(m).
  \notag
\end{equation}
\end{theorem}

A positive continuous function $u$ on $[0, \infty)$ is
called {\it (log, exp)-convex} if $\log u(e^{x})$ is
convex on $\spr$.

\begin{theorem} \label{thm:2-2}
Let $u\in C_{+, \log}$ be (log, exp)-convex. Then
\par\noindent
(1) $\ell_{u}(t)$ is decreasing for large $t$,
\par\noindent
(2) $\lim_{t\to\infty} \ell_{u}(t)^{1/t} =0$,
\par\noindent
(3) $u(r) = \sup_{t\geq 0} \ell_{u}(t) r^{t}$ for all
$r\geq 0$.
\end{theorem}

Let $k>0$. A positive continuous function $u$ on
$[0, \infty)$ is called {\it (log, $x^{k}$)-convex} if
$\log u(x^{k})$ is convex on $[0, \infty)$.

\begin{theorem} \label{thm:2-3}
Let $u\in C_{+, \log}$. We have the assertions:
\par\noindent
(1) $u$ is (log, $x^{k}$)-convex if and only if
$\ell_{u}(t) t^{kt}$ is log-convex.
\par\noindent
(2) If $u$ is (log, $x^{k}$)-convex, then for any
nonnegative integers $n$ and $m$,
\begin{equation}
\ell_{u}(n) \ell_{u}(m) \leq \ell_{u}(0) 2^{k(n+m)}
  \ell_{u}(n+m).  \notag
\end{equation}
\end{theorem}

Suppose $u\in C_{+, \log}$ satisfies the condition
$\lim_{n\to\infty} \ell_{u}(n)^{1/n}=0$. We define its
{\it $L$-function} by
\begin{equation} \label{eq:2-3}
\cl_{u}(r) = \sum_{n=0}^{\infty} \ell_{u}(n) r^{n}.
\end{equation}
Note that the function $\cl_{u}$ is an entire function.
By Theorem \ref{thm:2-2}(2), $\cl_{u}$ is defined for
any (log, exp)-convex function $u$ in $C_{+, \log}$.

\begin{theorem} \label{thm:2-4}
(1) Let $u\in C_{+, \log}$ be (log, exp)-convex. Then
its $L$-function $\cl_{u}$ is also (log, exp)-convex and
for any $a>1$,
\begin{equation}
\cl_{u}(r) \leq {ea \over \log a} u(ar), \qquad
\forall r\geq 0.  \notag
\end{equation}
\par\noindent
(2) Let $u\in C_{+, \log}$ be increasing and
(log, $x^{k}$)-convex. Then there exists a constant $C$,
independent of $k$, such that
\begin{equation}
u(r) \leq C \cl_{u}(2^{k} r), \qquad \forall r\geq 0.
    \notag
\end{equation}
\end{theorem}

Two positive functions $u$ and $v$ on $[0, \infty)$ are
called {\it equivalent} if there exist positive constants
$c_{1}, c_{2}, a_{1}, a_{2}$ such that
\begin{equation}
c_{1} u(a_{1}r) \leq v(r) \leq c_{2} u(a_{2}r), \qquad
\forall r\geq 0.  \notag
\end{equation}

It can be easily checked that if $u \in C_{+, \log}$ is
increasing and (log, $x^{k}$)-convex, then $u$ is
(log, exp)-convex. Thus by Theorem \ref{thm:2-4} the
functions $u$ and $\cl_{u}$ are equivalent for any
increasing (log, $x^{k}$)-convex function
$u\in C_{+, \log}$.

\medskip
Next, let $C_{+, 1/2}$ denote the set of all positive
continuous functions $u$ on $[0, \infty)$ satisfying the
condition:
\begin{equation}
\lim_{r\to\infty} {\log u(r) \over \sqrt{r}}=\infty. \notag
\end{equation}

The {\it dual Legendre transform} $u^{*}$ of $u \in
C_{+, 1/2}$ is defined to be the function
\begin{equation}
u^{*}(r) = \sup_{s\geq 0} {e^{2\sqrt{rs}} \over u(s)},
   \qquad  r\in [0, \infty). \notag
\end{equation}

Note that $C_{+, 1/2}\subset C_{+, \log}$. Below are some
properties of the dual Legendre transform. For the proofs,
see the paper \cite{akk3}.

\begin{theorem} \label{thm:2-5}
Let $u\in C_{+, 1/2}$. Then $u^{*}$ belongs to $C_{+, 1/2}$
and is an increasing (log, $x^{2}$)-convex function on
$[0, \infty)$.
\end{theorem}

\begin{theorem} \label{thm:2-6}
If $u\in C_{+, 1/2}$ is (log, $x^{2}$)-convex, then the
Legendre transform $\ell_{u^{*}}$ of $u^{*}$ is given by
\begin{equation}
\ell_{u^{*}} (t) = {e^{2t} \over \ell_{u}(t) t^{2t}},
\qquad t\geq 0.   \notag
\end{equation}
\end{theorem}

Suppose $u\in C_{+, 1/2}$ satisfies the condition
$\lim_{n\to\infty} \big(\ell_{u}(n) (n!)^{2}\big)^{-1/n}=0$.
We define its {\it $L^{\#}$-function} by
\begin{equation} \label{eq:2-7}
\cl_{u}^{\#}(r) = \sum_{n=0}^{\infty} {1 \over
  \ell_{u}(n)(n!)^{2}}\, r^{n}.
\end{equation}
Note that the function $\cl_{u}^{\#}$ is an entire function.
It can be checked that $\cl_{u}^{\#}$ is defined for any
(log, $x^{2}$)-convex function $u$ in $C_{+, 1/2}$.

\begin{theorem} \label{thm:2-7}
Let $u \in C_{+, 1/2}$ be (log, $x^{2}$)-convex. Then the
functions $u^{*}, \>\cl_{u^{*}}$, and $\cl_{u}^{\#}$ are
equivalent.
\end{theorem}

\section{CKS-space associated with a growth function}
  \label{sec:3}
\smallskip

Let $u\in C_{+, \log}$ be a fixed function. Define a
sequence of positive numbers by
\begin{equation} \label{eq:3-1}
\a (n) = {1 \over \ell_{u}(n) n!}, \qquad n\geq 0.
\end{equation}
Then from Equations (\ref{eq:1-1}) and (\ref{eq:2-3})
we see that
\begin{equation} \label{eq:3-a}
G_{1/\a}(r) = \sum_{n=0}^{\infty} {1 \over n!\a(n)} r^{n}
= \sum_{n=0}^{\infty} \ell_{u}(n) r^{n}=\cl_{u}(r).
\end{equation}
Moreover, by Equations (\ref{eq:1-1}) and (\ref{eq:2-7}),
we have
\begin{equation} \label{eq:3-b}
G_{\a}(r) = \sum_{n=0}^{\infty} {\a(n) \over n!} r^{n}
= \sum_{n=0}^{\infty} {1\over \ell_{u}(n) (n!)^{2}} r^{n}
= \cl_{u}^{\#} (r).
\end{equation}

Under various conditions on $u$, we will show that the
corresponding sequence $\{\a(n)\}$ satisfies the six
essential conditions (A1) (A2) (B2) ($\wt{\text{B}}2$)
(C2) (C3) stated in Section \ref{sec:1}. A weaker
condition than (B2), called near-(B2), will be defined
below.

\begin{lemma}
Let $u\in C_{+, \log}$. Suppose $\inf_{r>0} u(r)=1$ and
$u$ satisfies the condition
\begin{equation} \label{eq:a}
\lim_{r\to\infty} {\log u(r) \over r} < \infty.
\end{equation}
Then the sequence $\{\a(n)\}$ satisfies condition (A1).
\end{lemma}

\begin{pf}
From the definition of the Legendre transform and the
assumption, we get $\ell_{u}(0) =\inf_{r>0} u(r) =1$. Hence
$\a(0)=1$. It is easy to check that Equation (\ref{eq:a})
holds if and only if there exist constants $a, c >0$ such
that $u(r) \leq ce^{ar}$. Therefore,
\begin{equation}
\ell_{u}(n) = \inf_{r>0} {u(r) \over r^{n}} \leq
  \inf_{r>0} {ce^{ar} \over r^{n}} = ca^{n} \left(
  {e\over n}\right)^{n}.  \notag
\end{equation}
This implies that $\a(n) = (\ell_{u}(n) n!)^{-1} \geq
(n! c a^{n})^{-1} (n/e)^{n}$. But from the Stirling
formula, we have $n! \leq e2^{n/2} (n/e)^{n}$. Hence
$\a(n) \geq \big(ce(a\sqrt{2})^{n}\big)^{-1}$. We can
choose $\s=a\sqrt{2}$ to see that condition (A1) is
satisfied.
\end{pf}

\begin{lemma}
Suppose $u\in C_{+, 1/2}$ is (log, $x^{2}$)-convex. Then
the sequence $\{\a(n)\}$ satisfies condition (A2).
\end{lemma}

\begin{pf}
By using the Stirling formula we get
\begin{equation} \label{eq:3-2}
{\a(n) \over n!} = {1 \over \ell_{u}(n) (n!)^{2}}
\approx {1 \over 2\p n}{e^{2n} \over
\ell_{u}(n) n^{2n}}, \qquad \text{~for large~} n.
\end{equation}
By Theorem \ref{thm:2-5}, the dual Legendre transform
$u^{*}$ of $u$ belongs to $C_{+, 1/2}$ (hence to
$C_{+, \log}$) and is increasing and (log, $x^{2}$)-convex
(hence (log, exp)-convex.) Thus we can apply Theorem
\ref{thm:2-2} to $u^{*}$ and by part (2) of that theorem,
\begin{equation} \label{eq:3-3}
\lim_{n\to\infty} \ell_{u^{*}} (n)^{1/n} = 0.
\end{equation}
On the other hand, by Theorem 2.6
\begin{equation} \label{eq:3-4}
\ell_{u^{*}} (n) = {e^{2n} \over \ell_{u} (n) n^{2n}}.
\end{equation}
Therefore, from Equations (\ref{eq:3-2}) (\ref{eq:3-3})
(\ref{eq:3-4}), we easily see that $\big(\a(n)/n!
\big)^{1/n} \to 0$ as $n\to\infty$. Hence condition (A2)
is satisfied.
\end{pf}

Two sequences $\{a(n)\}$ and $\{b(n)\}$ are called
{\it equivalent} if there exist positive constants
$K_{1}, K_{2}, c_{1}, c_{2}$ such that for all $n$
\begin{equation}
K_{1}c_{1}^{n} a(n) \leq b(n)
   \leq  K_{2}c_{2}^{n} a(n).  \notag
\end{equation}

A sequence is said to be {\it nearly log-concave} if it
is equivalent to a log-concave sequence. With this
concept we define a weaker condition than (B2).

\medskip\noindent
$\bullet$ Near-(B2): The sequence $\{\a(n)/n!\}$ is nearly
log-concave.
\medskip

Observe that near-(B2) condition can be stated in another
way: The sequence $\{\a(n)\}$ is equivalent to a sequence
$\{\d(n)\}$ such that $\{\d(n)/n!\}$ is log-concave.

In Section \ref{sec:1}, we pointed out that condition
(B2) implies condition (B1). It can be easily checked
that condition near-(B2) also implies condition (B1).

\begin{lemma}
(1) Let $u\in C_{+, \log}$ be (log, $x^{1}$)-convex.
Then the sequence $\{\a(n)\}$ satisfies condition (B2).
\par\noindent
(2) Let $u\in C_{+, \log}$ be (log, $x^{2}$)-convex.
Then the sequence $\{\a(n)\}$ satisfies condition
near-(B2).
\end{lemma}

\begin{pf}
Let $u\in C_{+, \log}$ be (log, $x^{1}$)-convex. By
Theorem \ref{thm:2-3}(1) the sequence
$\{\ell_{u}(n)n^{n}\}$ is log-convex and so
\begin{equation}
\ell_{u}(n) n^{n} \ell_{u}(n+2) (n+2)^{n+2} \geq
\left(\ell_{u}(n+1) (n+1)^{n+1}\right)^{2}. \notag
\end{equation}
This inequality can be rewritten as
\begin{equation}
\ell_{u}(n+1)^{2} \leq {n^{n}(n+2)^{n+2} \over
(n+1)^{2(n+1)}} \ell_{u}(n) \ell_{u}(n+2). \notag
\end{equation}
But it can be easily verified that
\begin{equation}
{n^{n}(n+2)^{n+2} \over (n+1)^{2(n+1)}} \leq
\left({n+2 \over n+1}\right)^{2}. \notag
\end{equation}
Therefore,
\begin{equation}
\ell_{u}(n+1)^{2} \leq \left({n+2 \over n+1}\right)^{2}
\ell_{u}(n) \ell_{u}(n+2). \notag
\end{equation}
Rewrite this inequality as
\begin{equation}
{1 \over \ell_{u}(n)(n!)^{2}} {1 \over \ell_{u}(n+2)
((n+2)!)^{2}} \leq \left({1 \over
  \ell_{u}(n+1)((n+1)!)^{2}} \right)^{2}.
\end{equation}
Note that $\a(n)/n!=\big(\ell_{u}(n)(n!)^{2}\big)^{-1}$.
Hence the last inequality shows that the sequence
$\{\a(n)/n!\}$ is log-concave, i.e., the sequence
$\{\a(n)\}$ satisfies condition (B2). This proves the
first assertion of the lemma.

To prove the second assertion, let $u\in C_{+, \log}$
be (log, $x^{2}$)-convex. By Theorem \ref{thm:2-3}(1)
the sequence $\{\ell_{u}(n)n^{2n}\}$ is log-convex. But
$\a(n)=(\ell_{u}(n) n!)^{-1}$. Hence the sequence
$\{(\a (n) n!)^{-1} n^{2n}\}$ is log-convex and so the
sequence $\{\a(n)n!/n^{2n}\}$ is log-concave.

On the other hand, note that by the Stirling formula the
sequences $\{n!\}$ and $\{n^{n}\}$ are equivalent. This
implies that the sequence $\{\a(n)/n!\}$ is equivalent
to $\{\a(n)n!/n^{2n}\}$, which has just been shown to be
log-concave. Thus $\{\a(n)/n!\}$ is nearly log-concave
and so $\{\a(n)\}$ satisfies condition near-(B2).
\end{pf}

\begin{lemma}
Let $u\in C_{+, \log}$. Then the sequence $\{\a(n)\}$
satisfies condition ($\wt{\text{B}}2$).
\end{lemma}

\begin{pf}
By Theorem \ref{thm:2-1}(1) the sequence $\{\ell_{u}(n)\}$
is log-concave. Since $\ell_{u}(n)=(n!\a(n))^{-1}$, the
sequence $\{(n!\a(n))^{-1}\}$ is log-concave. Hence the
sequence $\{\a(n)\}$ satisfies condition ($\wt{\text{B}}2$).
\end{pf}

\begin{lemma}
Let $u\in C_{+, \log}$ be (log, $x^{2}$)-convex and
$\inf_{r>0} u(r) =1$. Then the sequence $\{\a(n)\}$
satisfies condition (C2).
\end{lemma}

\begin{pf}
From Theorem \ref{thm:2-3}(2) with $k=2$, we have
\begin{equation}
\ell_{u}(n) \ell_{u}(m) \leq \ell_{u}(0) 2^{2(n+m)}
  \ell_{u}(n+m).  \notag
\end{equation}
But $\ell_{u}(0)=\inf_{r>0} u(r)=1$ and $\a(n)=(\ell_{u}(n)
n!)^{-1}$. Hence
\begin{equation}
{1\over \a(n)n!} {1\over \a(m)m!} \leq 2^{2(n+m)}
  {1\over \a(n+m)(n+m)!}.  \notag
\end{equation}
Note that $n!m!/(n+m)! \leq 1$. Hence $\a(n+m) \leq
2^{2(n+m)} \a(n)\a(m)$ and so the sequence $\{\a(n)\}$
satisfies condition (C2).
\end{pf}

\begin{lemma}
Let $u\in C_{+, \log}$ and suppose $\inf_{r>0} u(r) =1$.
Then the sequence $\{\a(n)\}$ satisfies condition (C3).
\end{lemma}

\begin{pf}
By Theorem \ref{thm:2-1}(2) we have
\begin{equation}
\ell_{u}(0) \ell_{u}(n+m) \leq \ell_{u}(n) \ell_{u}(m).
  \notag
\end{equation}
Note that $\ell_{u}(n)=(\a(n) n!)^{-1}$ and $\ell_{u}(0)
= \inf_{r>0} u(r) =1$. Hence
\begin{equation}
\a(n) \a(m) \leq {(n+m)! \over n! m!} \a(n+m) \leq
2^{n+m} \a (n+m). \notag
\end{equation}
This shows that the sequence $\{\a(n)\}$ satisfies
condition (C3).
\end{pf}

Now, we summarize the above lemmas to state a theorem
which can be used for white noise distribution theory
on a CKS-space arising from a growth function.

\begin{theorem}
Let $u \in C_{+, 1/2}$ be a function satisfying the
following conditions:
\begin{itemize}
\item[(U0)] $\inf_{r>0} u(r) = 1$.
\smallskip
\item[(U2)] $\lim_{r\to\infty} r^{-1} \log u(r) < \infty$.
\smallskip
\item[(U3)] $u$ is (log, $x^{2}$)-convex.
\end{itemize}
\smallskip\noindent
Then the corresponding sequence $\{\a(n)\}$
satisfies conditions (A1), (A2), near-(B2),
($\wt{\text{B}}2$), (C2), and (C3).
\end{theorem}

\noindent
{\bf Remark}. Define another condition on $u$ by
\smallskip
\begin{itemize}
\item[(U1)] $u$ is increasing and $u(0)=1$.
\end{itemize}
\smallskip\noindent
Obviously, condition (U1) implies condition (U0). We need
this stronger condition in order to apply Theorem
\ref{thm:2-4}(2) where $u$ is assumed to be increasing.

\medskip
Now, let $u \in C_{+, 1/2}$ be a function satisfying
conditions (U0) (U2) (U3). From this function $u$ we define
a sequence $\{\a(n)\}$ of positive numbers by Equation
(\ref{eq:3-1}). Then we use this sequence $\{\a(n)\}$ to
define a CKS-space as in Equation (\ref{eq:b}). The
resulting Gel'fand triple is denoted by
\begin{equation}
[\ce]_{u} \subset (L^{2}) \subset [\ce]_{u}^{*}.
\end{equation}

Let us consider the characterization of generalized
functions in $[\ce]_{u}^{*}$. Note that by Equation
(\ref{eq:3-b}) the exponential generating function
$G_{\a}$ of $\{\a(n)\}$ is given by $G_{\a}=\cl_{u}^{\#}$.
Hence the growth condition in Equation (\ref{eq:1-5}) can
be stated as: There exist constants $K, a, p\geq 0$ such
that
\begin{equation}
|F(\x)| \leq K \cl_{u}^{\#} \big(a|\x|_{p}^{2}
  \big)^{1/2},  \qquad  \x\in\ce_{c}.  \notag
\end{equation}
Recall that by Theorem \ref{thm:2-7} the function
$\cl_{u}^{\#}$ is equivalent to the dual Legendre transform
$u^{*}$ of $u$. Thus we can replace the growth function
$\cl_{u}^{\#}$ by $u^{*}$ and we have the following theorem.

\begin{theorem}
Suppose $u\in C_{+, 1/2}$ satisfies the conditions (U0)
(U2) (U3). Then a complex-valued function $F$ on $\ce_{c}$
is the $S$-transform of a generalized function in
$[\ce]_{u}^{*}$ if and only if it satisfies the conditions:

\smallskip
\begin{itemize}
\item[(a)] For any $\x, \y\in \ce_{c}$, the function
$F(z\x+\y)$ is an entire function of $z\in\spc$.
\item[(b)] There exist constants $K, a, p\geq 0$, such
that
\begin{equation}
|F(\x)| \leq K u^{*}\big(a|\x|_{p}^{2}\big)^{1/2},
\qquad \x\in \ce_{c}.
\end{equation}
\end{itemize}
\end{theorem}

\medskip
Next, we consider the characterization of test functions in
$[\ce]_{u}$. By Equation (\ref{eq:3-a}) we have $G_{1/\a}=
\cl_{u}$. Hence the growth condition in Equation
(\ref{eq:1-6}) can be stated as: For any $a, p\geq 0$,
there exists a constant $K\geq 0$ such that
\begin{equation}
|F(\x)| \leq K \cl_{u} \big(a|\x|_{-p}^{2}\big)^{1/2},
  \qquad  \x\in\ce_{c}. \notag
\end{equation}

Assume that $u$ is increasing (so we need condition (U1).)
Then we can apply Theorem \ref{thm:2-4} to see that
$\cl_{u}$ is equivalent to the function $u$. Thus we can
replace the growth function $\cl_{u}$ by $u$ and we have
the next theorem.

\begin{theorem}
Suppose $u\in C_{+, 1/2}$ satisfies the conditions (U1)
(U2) (U3). Then a complex-valued function $F$ on $\ce_{c}$
is the $S$-transform of a test function in $[\ce]_{u}$
if and only if it satisfies the conditions:

\smallskip
\begin{itemize}
\item[(a)] For any $\x, \y\in \ce_{c}$, the function
$F(z\x+\y)$ is an entire function of $z\in\spc$.
\item[(b)] For any constants $a, p\geq 0$, there exists a
constant $K\geq 0$ such that
\begin{equation}
|F(\x)| \leq K u\big(a|\x|_{-p}^{2}\big)^{1/2},
\qquad \x\in \ce_{c}.
\end{equation}
\end{itemize}
\end{theorem}

\medskip
At the end of this paper we give two examples to illustrate
our method.

\begin{example}
Consider the Kondratiev-Streit space \cite{ks92} \cite{ks93}
\cite{kuo96} associated with the sequence $\a(n)=(n!)^{\b},
0\leq \b<1$. If we apply the characterization theorems from
\cite{cks} and \cite{akk2}, then we have to use the
following growth functions:
\begin{equation} \label{eq:3-5}
G_{\a}(r) = \sum_{n=0}^{\infty} {1\over (n!)^{1-\b}}\,
r^{n}, \qquad  G_{1/\a}(r) = \sum_{n=0}^{\infty}
{1\over (n!)^{1+\b}}\, r^{n}.
\end{equation}
But these growth functions are impractical since the series
cannot be summed up in closed forms when $\b\not= 0$. To
overcome this difficulty, note that the sequence $\{n!\}$ is
equivalent to the sequence $\{(n/e)^{n}\}$ by the Stirling
formula. Hence the function $G_{1/\a}$ is equivalent to the
following function
\begin{equation}
  \t (r) = \sum_{n=0}^{\infty} \left({e\over n}
   \right)^{(1+\b)n} r^{n}. \notag
\end{equation}
As shown in \cite{akk3}, $(e/n)^{(1+\b)n} = \ell_{u}(n)$
for the function
\begin{equation}
  u(r)=\exp\left[(1+\b) r^{1/(1+\b)}\right]. \notag
\end{equation}
Therefore $\t=\cl_{u}$ by Equation (\ref{eq:2-3}) and so
$G_{1/\a}$ is equivalent to $\cl_{u}$. Furthermore, by
Theorem \ref{thm:2-4}, the function $\cl_{u}$ is equivalent
to $u$. Hence the growth function $G_{1/\a}$ in Equation
(\ref{eq:3-5}) can be replaced by the above function $u$.

We remark that the Kondratiev-Streit space turns out to be
the same as the CKS-space arising from the function $u$.
It is easy to check that the dual Legendre transform
$u^{*}$ of $u$ is given by
\begin{equation}
u^{*} (r)=\exp\left[(1-\b) r^{1 \over (1-\b)}\right].
  \notag
\end{equation}
Hence the growth function $G_{\a}$ in Equation
(\ref{eq:3-5}) can be replaced by $u^{*}$. Thus we have
derived the growth functions $u$ and $u^{*}$ used by
Kondratiev and Streit \cite{ks92} \cite{ks93}.
\end{example}

\begin{example}
Let $u(r) = \exp_{k}(r) = \exp(\exp(\cdots (\exp(r))))$
be the $k$-th iteration of the exponential function. It
is shown in \cite{akk3} that the dual Legendre transform
$u^{*}$ of $u$ is equivalent to the following function
\begin{equation}
v(r) = \exp\left[2 \sqrt{r \log_{k-1}\sqrt{r}}\,\right],
   \notag
\end{equation}
where $\log_{j}$ is the function defined by
\begin{equation}
\log_{1}(r)=\log(\max\{r, e\}), \quad  \log_{j}(r)=
\log_{1}(\log_{j-1}(r)), \quad j\geq 2.  \notag
\end{equation}
The CKS-space arising from the function $v$ turns out to
be the same as the one defined by the Bell numbers
$\{b_{k}(n)\}$ of order $k$. The function $v$ serves as a
growth function for the characterization theorem of test
functions in this CKS-space. The function $v^{*}$ is
equivalent to the function $u(r) = \exp_{k}(r)$. Hence
$\exp_{k}$ can be used as a growth function for the
characterization theorem of generalized functions in this
CKS-space.
\end{example}

\bigskip
\noindent
{\bf Acknowledgements.} N.~Asai wants to thank the Daiko
Foundation and the Kamiyama Foundation for research support.
H.-H.~Kuo is grateful for financial supports from the
Academic Frontier in Science (AFS) of Meijo University.
He wants to give his deepest appreciation to Professors
T. Hida and K. Sait\^o for the warm hospitality during his
visits (February 15--21, 1998 and March 1--7, 1999.)

\end{document}